\documentclass[12pt,twoside,leqno]{amsart}
\usepackage{mathptmx}
\usepackage{amssymb, amsmath, amsthm, enumerate}
\usepackage{graphicx}
\usepackage{tikz}
\begin{document}
\setcounter{page}{1}
\vspace*{2.0cm}

\title[ON THE SOMBOR INDEX AND SOMBOR ENERGY OF $m$-SPLITTING GRAPH AND $m$-SHADOW GRAPH OF REGULAR GRAPHS] {ON THE SOMBOR INDEX AND SOMBOR ENERGY OF $m$-SPLITTING GRAPH AND $m$-SHADOW GRAPH OF REGULAR GRAPHS}
\author{Randhir Singh$^*$}
\author {S.C. Patekar }
\date{}
\maketitle

\vspace*{-0.3in}

\begin{center}
{\footnotesize  Department of Mathematics\\ Savitribai Phule Pune University, Pune 411007, India\\
\email{randhir\textunderscore singh46@yahoo.com}, \email{shri82patekar@gmail.com}
}

\end{center}
\vskip 2mm

{\footnotesize 

\begin{abstract}
	
 A vertex-degree-based topological index named as Sombor index of a simple graph $G$ with $n$ vertices was recently introduced by I. Gutman \cite{IG sombor index} and defined as $SO=SO(G)=\sum\limits_{{uv}\in E(G)}\sqrt{{d_{G}(u)}^2+{d_{G}(v)}^2}$, where $d_{G}(u)$ is degree of vertex $u$ and $E(G)$ is edge set. In this paper, we find Sombor index of $m$-splitting graph and $m$-shadow graph. Also, we determine relation between energy and
 Sombor energy of $ m $-splitting graph and $ m $-shadow graph of $k$-regular graph. \vskip 1mm
\end{abstract}
\noindent {\bf Keywords:} $m$-Splitting Graph, $m$-Shadow Graph, Sombor index, Sombor energy.\vskip 1mm

%---- http://www.ams.org/msc/msc2010.html ---%
\noindent {\bf 2010 AMS Subject Classification:}  05C50, 05C09

\vskip 6mm

\vskip 6mm
\noindent\textbf{1. Introduction}
\vskip 6mm

\noindent Let $G$ be a simple graph with vertex set $V(G)$ and edge set $E(G)$, $|V(G)|=n$. If the vertices $u,v\in V(G)$ are adjacent, then edge whose end points are $u$ and $v$ is denoted by $uv$. The neighbor of $u$ is denoted by $N(u)$ that is the set of vertices adjacent to $u$, $|N(u)|$ is total number of adjacent vertices to $u$ and is called degree of $u$ and denoted by $d_{G}(u)$.

Let $A(G)$ be the adjacency matrix of a simple graph $G$ with vertices $v_1,v_2,...,v_n$, elements of adjacency matrix are defined by
\begin{center}
		$A(G) = (a_{ij})=  \begin{cases}
		$1$ & \text{if $v_i$ and $v_j$ are adjacent}\\
		$0$ & \text {if $i=j$}\\
		$0$ & \text {if $v_i$ and $v_j$ are not adjacent}.\\
	\end{cases}$
\end{center}
and let $\lambda_1, \lambda_2,...,\lambda_n$ be the eigenvalues of matrix $A(G)$.

The energy of simple graph, introduced by I. Gutman \cite{IG EG}, is defined as
\begin{center}
	$\varepsilon(G)=\sum\limits_{i=1}^{n}|\lambda_i|$
\end{center}
For graph-theoretical notions we refer \cite{Babat,DBW}.

\vspace{0.1cm}
\noindent\textbf{1. Sombor index of $m-$ Splitting Graph.}
\vspace{0.1cm}

\noindent\textbf{Definition 1.1} \cite{sampathkumar}
The \textit{splitting graph} $S'(G)$ of a graph $G$ is obtained by adding to each vertex $v$ a new vertex $v'$ such that $v'$ is adjacent to every vertex that is adjacent to $v$ in $G$.

\noindent\textbf{Definition 1.2} 
The \textit{$m$-splitting graph} $Spl_{m}(G)$ of a graph $G$ is obtained
by adding to each vertex $v_i$ of $G$ new $m$ vertices, say $v_{ij}$ such that $v_{ij}$, $1 \leqslant i\leqslant n$ and $1 \leqslant j\leqslant  m$, is adjacent to each vertex that is adjacent to $v_i$ in $G$.

Now we give general formula for the Sombor index of $k$-regular graph and $m$-splitting graph of $k$-regular graph.

\noindent\textbf{Theorem 1.}\textit{ Let $G$ is a $k$-regular graph with $n$, where $n\geqslant2$, vertices and $Spl_{m}(G)$ is $m$-splitting graph of $G$. Then}
\begin{enumerate}
	\item \textit{$SO(G)=\frac{nk^2}{\sqrt{2}}$,}  
	\item \textit{ $SO(Spl_{m}(G))=SO(G)(m\sqrt{2m^2+4m+4}+(m+1))$.}
\end{enumerate}
\noindent\textit{Proof.} \begin{enumerate}
	\item 
 Let $G$ be $k$-regular graph with $n$ vertices.

$\Rightarrow e(G)=\frac{nk}{2}$

$\because SO(G)=\sum\limits_{{v_{i}v_{j}}\in E(G)}\sqrt{d_{G}({v_i})^{2}+{d_{G}({v_{j}})^{2}}}$ \space from \cite{IG sombor index}

and $d_{v_i}=k $ \space \space for all $v_i\in V(G)$.

$\therefore SO(G)=e(G)\sqrt{k^2+k^2}$

\hspace{0.64in}$=\frac{nk}{2}\sqrt{2k^2}$

\hspace{0.64in}$=\frac{nk^2}{\sqrt{2}}$.

\item $|V(Spl_{m}(G))\backslash V(G)|=nm$ and $d_{G}(v_{ij})=k$ \space \space for all $1\leqslant i \leqslant n$ and $1\leqslant j\leqslant m$. $d_{G}(v_{i})$ for all $1\leqslant i\leqslant n$ increased by $mk$ therefore $d_{G}(v_{i})=(m+1)k$.

$\because N(v_{ij})=\{v_l:v_l \text{ is adjacent to } v_i\}$ \space \space for all $i\neq l$ and $1\leqslant j\leqslant m$ and $d_{Spl_{m}(G)}(v_{l})=(m+1)k$.

$\therefore SO(Spl_{m}(G))=\sum\limits_{v_{ij}v_{l}\in E(Spl_{m}(G))}\sqrt{(d_{Spl_{m}(G)}(v_{ij}))^2+(d_{Spl_{m}(G)}(v_{l}))^2}$\\

\vspace{-0.5cm}

 \hspace{2.5 in} $+\sum\limits_{v_{i}v_{j}\in E(Spl_{m}(G))}\sqrt{(d_{Spl_{m}(G)}(v_{i}))^2+(d_{Spl_{m}(G)}(v_{j}))^2}$

$=|N(v_{ij})||V(Spl_{m}(G))\backslash V(G)|\sqrt{k^2+((m+1)k)^2}+e(G)\sqrt{((m+1)k)^2+((m+1)k)^2}$

$=knm\sqrt{k^2(m^2+2m+2)} + \frac{n(m+1)k^2}{2}\sqrt{2}$

$=\frac{nk^2}{\sqrt{2}}[m\sqrt{2m^2+4m+4}+(m+1)]$

$=SO(G)[m\sqrt{2m^2+4m+4}+(m+1)]$.
\end{enumerate}

\noindent\textbf{Illustration.} Consider complete graph $K_n$ which is $(n-1)$-regular.

$SO(Spl_m(K_n))=\frac{n(n-1)^2}{\sqrt{2}}[m\sqrt{2m^2+4m+4}+(m+1)]$.

\noindent\textbf{Example.} Consider $3-$regular graph with $6$ vertices

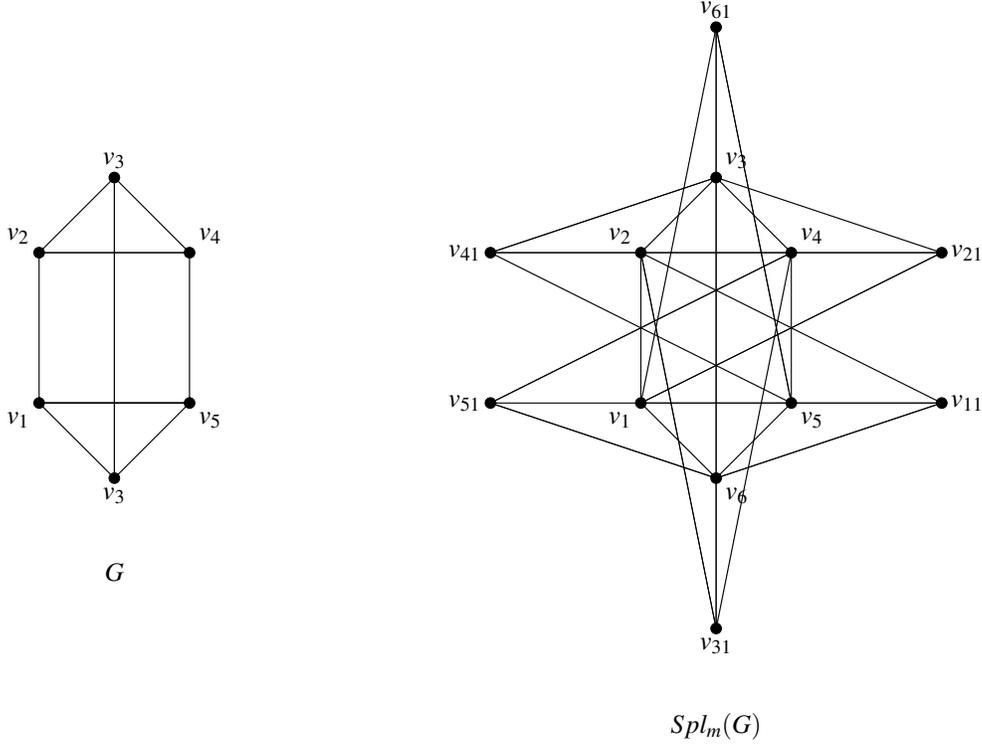
\begin{tikzpicture}
\filldraw [black] (0,0) circle (2pt) node[below left] {$v_1$}--(0,2) circle (2pt) node[above left] {$v_2$}--(1,3) circle (2pt) node[above]{$v_3$}--(2,2)circle (2pt) node[ above right]{$v_4$}--(2,0)circle (2pt) node[below right]{$v_5$}--(1,-1)circle (2pt) node[below]{$v_3$} --(0,0);
\filldraw (0,0)--(2,0);
\filldraw (1,-1)--(1,3);
\filldraw (2,2)--(0,2);
\filldraw (0,0)--(2,0);
\filldraw[black] (8,0) circle (2pt) node [below left] {$v_1$}--(9,-1) circle (2pt) node [below right]{$v_6$}--(10,0) circle (2pt) node [below right]{$v_5$}--(10,2) circle (2pt) node [above right]{$v_4$}--(9,3) circle (2pt) node [above right]{$v_3$}--(8,2) circle (2pt) node [above left]{$v_2$}--(8,0)--cycle ;
\filldraw (8,0)--(10,0);
\filldraw (9,-1)--(9,3);
\filldraw (8,2)--(10,2);
\filldraw[black] (6,0) circle (2pt) node[left]{$v_{51}$}--(8,0);
\filldraw (6,0)--(9,-1);
\filldraw (6,0)--(10,2);
\filldraw[black] (6,2) circle (2pt) node [left]{$v_{41}$}--(10,0);
\filldraw (6,2)--(8,2);
\filldraw (6,2)--(9,3);
\filldraw[black] (9,5) circle (2pt) node
[above] {$v_{61}$}--(8,0);
\filldraw (9,5)--(9,3);
\filldraw (9,5)--(10,0);
\filldraw[black] (12,2) circle (2pt) node
[right] {$v_{21}$}--(9,3);
\filldraw (12,2)--(10,2);
\filldraw (12,2)--(8,0);
\filldraw[black] (12,0) circle (2pt) node
[right] {$v_{11}$}--(8,2);
\filldraw (12,0)--(10,0);
\filldraw (12,0)--(9,-1);
\filldraw[black] (9,-3) circle (2pt) node
[below] {$v_{31}$}--(10,2);
\filldraw (9,-3)--(9,-1);
\filldraw (9,-3)--(8,2);
\draw (1,-2) node[below]{$G$};
\draw (9,-4) node[below]{$Spl_m(G)$};

\end{tikzpicture}

Here $m=1$, $k=3$, and $n=6$.\\
$\displaystyle SO(G)=\sum\limits_{{v_{i}v_{j}}\in E(G)}\sqrt{d_{G}({v_i})^{2}+{d_{G}({v_{j}})^{2}}}$\\

\hspace{0.28in} $\displaystyle =9\sqrt{3^2+3^2}$\\

\hspace{0.28in} $\displaystyle =27\sqrt{2}$.\\

$\displaystyle SO(Spl_m(G))= \sum\limits_{{v_{ij}v_{l}}\in E(Spl_m(G))}\sqrt{d_{Spl_m(G)}({v_{ij}})^{2}+{d_{Spl_m(G)}({v_{l}})^{2}}}+\sum\limits_{{v_{i}v_{j}}\in E(Spl_m(G))}\sqrt{d_{Spl_m(G)}({v_i})^{2}+{d_{Spl_m(G)}({v_{j}})^{2}}}$\\

\hspace{0.9in} $= 18\sqrt{3^2+6^2}+9\sqrt{6^2+6^2 }$\\

\hspace{0.9in} $=54(\sqrt{5}+\sqrt{2})=SO(G)(\sqrt{10}+2)$.

\noindent\textbf{2. Sombor index of $m-$ Shadow Graph.}

\noindent\textbf{Definition 1.2.1.} The \textit{shadow graph} $D_2(G)$ of a connected graph $G$ is constructed by taking two copies of $G$, say $G'$ and $G''$. Join each vertex $u'$ in $G'$ to the neighbors of the corresponding vertex $u''$ in $G''$.

\noindent\textbf{Definition 1.2.2.} The \textit{$m$-shadow} graph $D_{m}(G)$ of a connected graph $G$ is constructed by taking $m$ copies of $G$, say $G_1, G_2, ..., G_m$ then join each vertex $u$ in $G_i$ to the neighbors of the corresponding vertex $v$ in, $G_j$, $1\leqslant i$, $j\leqslant m$.

\noindent\textbf{Preposition 1.} \textit{If $G$ is $k$-regular graph with $n$ vertices, then $m$-shadow graph $D_{m}(G)$ of $G$ is $mk$-regular.}

\noindent\textbf{Theorem 2.} If $G$ is $k$-regular graph with $n$ vertices, then $SO(D_{m}(G))= SO(G)(m^3+m^2)$.

\noindent\textit{Proof.} Since $D_{m}(G)$ is constructed from $m$ copies of $G$ and in $V(G)=n$. Therefore $|V(D_{m}(G))|=mn+n=(m+1)n$ and $e(G)=\frac{(nm+n)(mk)}{2}$.\\
$\therefore SO(D_{m}(G))= \sum\limits_{{uv}\in E(D_{m}(G))}\sqrt{{(d_{D_{m}(G)}(u))}^2+{(d_{D_{m}(G)}(v))}^2}$

\noindent Since sum runs over number of edges and degree of each end vertex of each edge is same and is equal to $mk$.
Therefore, 

$SO(D_{m}(G))= e(G)\sqrt{(mk)^2+(mk)^2}$

\hspace{0.81in}$= \frac{(nm+n)(mk)}{2}\sqrt{2(mk)^2}$

\hspace{0.81in}$=\frac{(nm+n)(mk)^2}{\sqrt{2}}$

\hspace{0.81in}$=\frac{nk^2}{\sqrt{2}}(m^3+m^2)$

\hspace{0.81in}$=SO(G)(m^3+m^2).$

\noindent\textbf{Illustration.} Consider cycle graph $C_n$ which is $2$-regular.

$SO(D_m(C_n))=\frac{4n}{\sqrt{2}}(m^3+m^2).$

\noindent\textbf{Example.} Consider cycle graph $C_6$.\\
Here, $m=1$ $n=6$ and $k=2.$\\
$SO(C_6)= 6\sqrt{2^2+2^2}$\\
$=12\sqrt{2}$\\
$SO(D_m(C_6))=\sum\limits_{{uv}\in E(D_{m}(G))}\sqrt{{(d_{D_{m}(G)}(u))}^2+{(d_{D_{m}(G)}(v))}^2}$\\
$= e(G)\sqrt{(mk)^2+(mk)^2}$\\
	$=12\sqrt{2^2+2^2}$\\
$=24\sqrt{2}=2SO(C_6)$.\\
\newpage
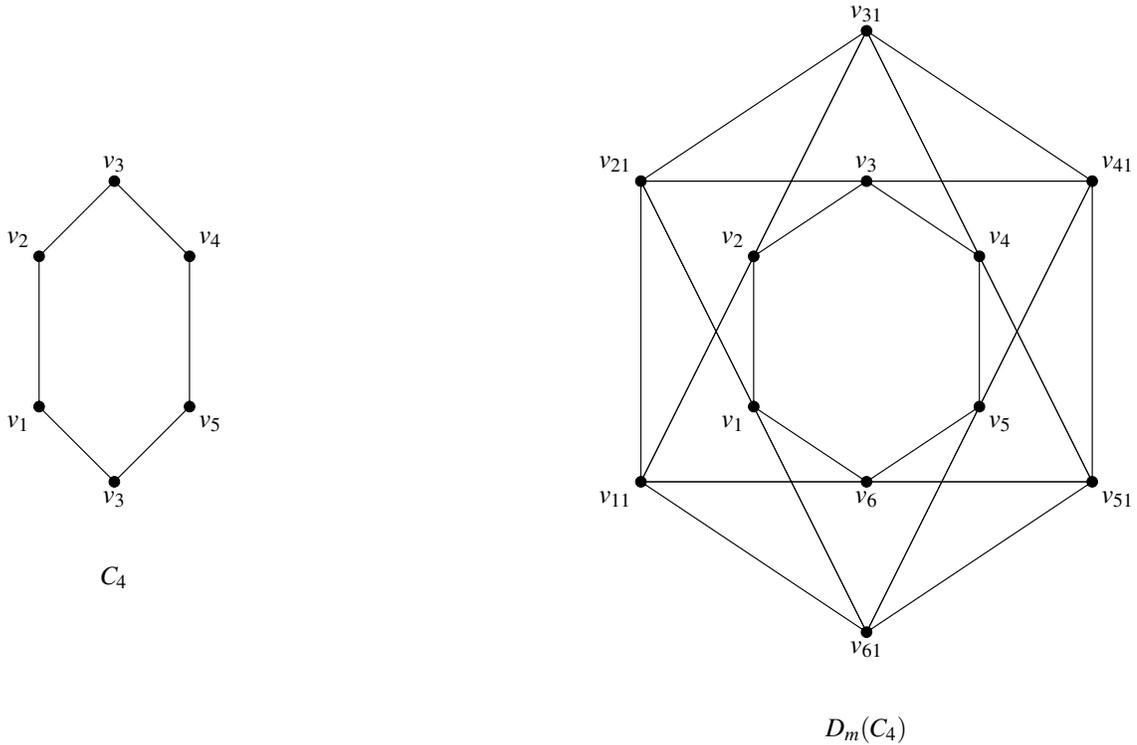
\begin{tikzpicture}
	\filldraw [black] (0,0) circle (2pt) node[below left] {$v_1$}--(0,2) circle (2pt) node[above left] {$v_2$}--(1,3) circle (2pt) node[above]{$v_3$}--(2,2)circle (2pt) node[ above right]{$v_4$}--(2,0)circle (2pt) node[below right]{$v_5$}--(1,-1)circle (2pt) node[below]{$v_3$} --(0,0);
	
	\filldraw [black] (9.5,0) circle (2pt) node[below left] {$v_1$}--(9.5,2) circle (2pt) node[above left] {$v_2$}--(11,3) circle (2pt) node[above]{$v_3$}--(12.5,2)circle (2pt) node[ above right]{$v_4$}--(12.5,0)circle (2pt) node[below right]{$v_5$}--(11,-1)circle (2pt) node[below]{$v_6$}--(9.5,0);

		\filldraw [black] (8,-1) circle (2pt) node[below left] {$v_{11}$}--(8,3) circle (2pt) node[above left] {$v_{21}$}--(11,5) circle (2pt) node[above]{$v_{31}$}--(14,3)circle (2pt) node[ above right]{$v_{41}$}--(14,-1)circle (2pt) node[below right]{$v_{51}$}--(11,-3)circle (2pt) node[below]{$v_{61}$}--(8,-1);
		\draw (1,-2) node[below]{$C_4$};
	\filldraw (8,-1)--(11,-1);
	\filldraw (8,-1)--(9.5,2);
	\filldraw (11,-3)--(9.5,0);
	\filldraw (11,-3)--(12.5,0);
	\filldraw (14,-1)--(11,-1);
	\filldraw (14,-1)--(12.5,2);
	\filldraw (14,3)--(12.5,0);
	\filldraw (14,3)--(11,3);
	\filldraw (11,5)--(9.5,2);
	\filldraw (11,5)--(12.5,2);
	\filldraw (8,3)--(9.5,0);
	\filldraw (8,3)--(11,3);
		\draw (11,-4) node[below]{$D_m(C_4)$};
	
\end{tikzpicture}

\noindent\textbf{3. Sombor Energy of $m$-splitting and $m$-shadow graph.}

Recently K. J. Gowtham et al. \cite{KJG} introduced a new matrix for a graph $G$ called the Sombor matrix and defined a new variant of graph energy called Sombor energy $ES(G)$ of graph $G$.

\noindent\textbf{Definition 3.1.} \cite{KJG} The Sombor matrix of a graph $G$ with $n$ vertices is defined as $S(G)=(s_{ij})_{n\times n},$ where 
\begin{center}
	$s_{ij}=$ $\begin{cases}
		\sqrt{{d_{G}(u)}^2+{d_{G}(v)}^2} &  \text{${uv}\in E(G)$}\\
		$0$ & $otherwise$\\
	\end{cases}$
\end{center}

Now we relate energy of regular graph $\varepsilon(G)$ with Sombor energy regular graph $ES(G)$ which is sum of absolute eigenvalues of Sombor matrix.

Let $A\in\mathbb{R}^{m\times n}$ and $B\in\mathbb{R}^{p\times q}$. Then the \textit{Kronecker Product} of $A$ and $B$ is defined as the matrix
\begin{center}
$A\otimes B=	\begin{bmatrix}
		a_{11}B & a_{12}B & ... & a_{1n}B\\
		a_{21}B & a_{22}B & ... & a_{2n}B\\
		\vdots & \vdots & \ddots & \vdots\\
		a_{m1}B & a_{m2}B & ... & a_{mn}B
	\end{bmatrix}_{mp\times nq}$
\end{center}

\noindent\textbf{Preposition 2.} \cite{horn}\textit{ Let $A \in M^m$ and $B\in Mn$. Furthermore, let $\lambda$ be an eigenvalues of matrix $A$ with corresponding eigenvector $x$, and $\mu$ an eigenvalue of matrix $B$ with corresponding eigenvector $y$. Then $\lambda \mu$ is an eigenvalue of $A \otimes B$ with corresponding eigenvector $x \otimes y$. }

\noindent\textbf{3.1. Sombor Energy of $m-$splitting graph.}  

\noindent\textbf{Theorem 3.} \textit{If $G$ is $k-$regular graph, then $SE(Spl_m(G))=k(m+1)\sqrt{2}\varepsilon(G)$}.

\noindent\textit{Proof.} Let $G$ be $k-$regular graph with vertex set $V(G)=\{v_1, v_2, v_3, ..., v_n\}$ and $Spl_m(G)$ be splitting graph of graph $G$ with vertex set $V(Spl_m(G))=\{v_{11}, v_{12},...,v_{1n}, v_{21}, v_{22},...,v_{2n},...v_{m1}, v_{m2},...,v_{mn}\}\cup V(G)$. Then Sombor matrix of $Spl_m(G)$ is
\begin{center}
	$S(Spl_m(G))=\begin{array}{ccccccccccc}&
		\begin{array}{c c c c ccccccc}
		~~	v_1~~ &~~ v_2~~ &~~ ... ~~ &~~ v_n ~~&~~ v_{11} &~...~&v_{1n} & ~... & v_{m1}& ...& v_{mn}\\
		\end{array}\\
		\begin{array}{ccccccccccc}
v_1 \\ v_2\\ \vdots\\ v_n \\ v_{11}\\ \vdots\\v_{1n} \\ \vdots \\ v_{m1} \\ \vdots\\ v_{mn}
		\end{array}
		&
		\left[
		\begin{array}{cccc|ccccccc}
			\begin{array}{cccc}
				0 & a_{11} & ... & a_{1n} \\
				a_{21} & 0 & ... & a_{2n}\\
				\vdots & \vdots&\ddots&\vdots\\
				a_{n1}& a_{n2}&...& 0\\				
	\hline
			0&b_{12}&...&b_{1n}\\
			\vdots&\vdots&...&\vdots\\
			b_{n1}&b_{n2}&...&0\\
			\vdots&\vdots&...&\vdots\\
			b_{m1}&b_{m2}&...&0\\
			\vdots&\vdots&...&\vdots\\
			b_{m1}&b_{m2}&...&0
			
		\end{array}
	&
	\vline
	&
	\begin{array}{ccccccc}
		0&...&b_{1n}&...&0&...&b_{1n}\\
		b_{21}&...&b_{2n}&...&b_{21}&...&b_{2n}\\
		\vdots&...&\vdots&...&\vdots&...&\vdots\\
		b_{n1}&...&0&...&b_{n1}&...&0\\
		\hline
		0&...&0&...&0&...&0\\
		\vdots&...&\vdots&...&\vdots&...&\vdots\\
		0&...&0&...&0&...&0\\
		\vdots&...&\vdots&...&\vdots&...&\vdots\\
		0&...&0&...&0&...&0\\
		\vdots&...&\vdots&...&\vdots&...&\vdots\\
		0&...&0&...&0&...&0
	\end{array}

		\end{array}
		\right]
	\end{array}$
\end{center}

which is $(mn+n)\times(mn+n)$ matrix. Sombor matrix $S(Spl_m(G))$ can be written as
\vspace{0.5cm}

\hspace{1.75cm}	$S(Spl_m(G))=\begin{bmatrix}	
	S_1&S_2&...&S_2\\
	S_2&0&...&0\\
	\vdots&\vdots&\ddots&\vdots\\
	S_2&0&...&0
\end{bmatrix}$

Where $S_1$ and $S_2$ are $n\times n$ symmetric matrices and entries of matrix $S_1$ are
\begin{center}
	$a_{ij}=\begin{cases}
		(m+1)k\sqrt{2} & \text{if $v_i$ is adjacent to $v_j$ }\\
		0 & \text{Otherwise}.
	\end{cases}$
\end{center}
Entries of $S_2$ are
\begin{center}
		$b_{ij}=\begin{cases}
		k\sqrt{m^2+2m+2} & \text{if $v_i$ is adjacent to $v_{ij}$ }\\
		0 & \text{Otherwise}.
	\end{cases}$

\end{center}

$\implies S_1= (m+1)k\sqrt{2}~A(G) \text{ and } S_2 = k\sqrt{m^2+2m+2}~A(G).$

\vspace{0.5 cm}

$\therefore$  \hspace{1.75cm} $S(Spl_m(G))=\begin{bmatrix}	
	(m+1)k\sqrt{2}~A(G)&k\sqrt{m^2+2m+2}~A(G)&...&k\sqrt{m^2+2m+2}~A(G)\\
	k\sqrt{m^2+2m+2}~A(G)&0&...&0\\
	\vdots&\vdots&\ddots&\vdots\\
	k\sqrt{m^2+2m+2}~A(G)&0&...&0
\end{bmatrix}_{(mn+n)}$

\hspace{1.68in} $=\begin{bmatrix}	
	(m+1)k\sqrt{2}&k\sqrt{m^2+2m+2}&...&k\sqrt{m^2+2m+2}\\
	k\sqrt{m^2+2m+2}&0&...&0\\
	\vdots&\vdots&\ddots&\vdots\\
	k\sqrt{m^2+2m+2}&0&...&0
\end{bmatrix}_{(m+1)} \otimes A(G)$

We are interested to find Sombor energy of splitting graph $Spl_m(G)$ therefore we need to find eigenvalues of matrix $S(Spl_m(G))$, for this
let $B= \begin{bmatrix}	
	(m+1)k\sqrt{2}&k\sqrt{m^2+2m+2}&...&k\sqrt{m^2+2m+2}\\
	k\sqrt{m^2+2m+2}&0&...&0\\
	\vdots&\vdots&\ddots&\vdots\\
	k\sqrt{m^2+2m+2}&0&...&0
\end{bmatrix}_{(m+1)}$. 
Let $\mu_1, \mu_2, ...,\mu_{m+1}$ be eigenvalues of matrix $B$. $B$ is of rank two therefore $A$ has two non-zero eigenvalues $\mu_1$ and $\mu_2$ (say).\\
For convenience set $a =(m+1)k\sqrt{2}$ and $b=k\sqrt{m^2+2m+2}$. Then\\

\begin{center}
		$tr(A)=\mu_1+\mu_2$\\
\begin{equation} \therefore~~~~ \mu_1+\mu_2=a \end{equation}
and\\

 $tr(A^2)=\mu_1^2+\mu_2^2$
\begin{equation}  
	\therefore~~~~\mu_1^2+\mu_2^2=2mb^2+a^2 \end{equation} 
\end{center}
Solving $ (1) $ and $ (2) $ we get\\
\begin{center}
$\mu^2 -a\mu-mb^2=0 $\\
\end{center}
whose roots are $ \mu_1 $ and $ \mu_2 $.

 $\therefore$ Characteristic equation of $B$ is\\
	$Char(A)=\mu^{m-1}((\mu^2 -k(m+1)\sqrt{2})\mu-mk^2(m^2+2m+2))=0$\\
	$\displaystyle \mu^{m-1}((\mu^2 -k(m+1)\sqrt{2})\mu-mk^2(m^2+2m+2))=0$\\
	$\displaystyle\implies \mu=0$\\ and
	$\displaystyle \mu= \frac{k(m+1)\sqrt{2}\pm k\sqrt{4m^3+10m^2+10m+2}}{2}$.\\
	Therefore,\begin{center}
		 $\displaystyle spec(B)=\begin{array}{ccc}
		
		\left(
		\begin{array}{ccc}
		0&\frac{k(m+1)\sqrt{2}+k\sqrt{4m^3+10m^2+10m+2}}{2}&\frac{k(m+1)\sqrt{2}- k\sqrt{4m^3+10m^2+10m+2}}{2}\\
		m-1& 1&1
	\end{array}
		\right)
	\end{array}$
	\end{center}

Therefore by preposition [2]
\begin{center}
	
	$\displaystyle spec(S(Spl_m(G)))=\begin{array}{cccccccccc}
        \left(
        \begin{array}{cccccccccc}
        	0\lambda_1&...&0\lambda_n&...&a\lambda_1&...&a\lambda_n&b\lambda_1&...&b\lambda_n\\
        	m-1&...&m-1&...&1&...&1&1&...&1
        \end{array}
        \right)
	\end{array}$
\end{center}

Where $\displaystyle a=\frac{k(m+1)\sqrt{2}+k\sqrt{4m^3+10m^2+10m+2}}{2}$

  and $\displaystyle
b=\frac{k(m+1)\sqrt{2}-k\sqrt{4m^3+10m^2+10m+2}}{2}$.\\
Therefore,

	$\displaystyle ES(Spl_m(G))=\sum\limits_{i=1}^{mn+n}\left|\frac{k(m+1)\sqrt{2}\pm k\sqrt{4m^3+10m^2+10m+2}}{2}\lambda_i\right|$\\
	
	$\displaystyle=\sum\limits_{i=1}^{n}\left|\lambda_i\right|\left[\frac{k(m+1)\sqrt{2}+ k\sqrt{4m^3+10m^2+10m+2}}{2}+\frac{k(m+1)\sqrt{2}- k\sqrt{4m^3+10m^2+10m+2}}{2}\right]$\\
$=\sum\limits_{i=1}^{n}|\lambda_i||k(m+1)\sqrt{2}|$\\
$=k(m+1)\sqrt{2}\varepsilon(G)$\\
\noindent\textbf{3.2. Sombor Energy of $m-$shadow graph.}
 
\noindent\textbf{Theorem 4.} If $G$ is $k-$regular graph. Then $ES(D_m(G))=mk\sqrt{2+8m}\varepsilon(G)$.\\
\noindent\textit{Proof.} Let $G$ be $k-$regular graph. Then Sombor matrix  is

$S(D_m(G))=\begin{bmatrix}
	mk\sqrt{2}A(G)&mk\sqrt{2}A(G)&...&mk\sqrt{2}A(G)\\
	mk\sqrt{2}A(G)&0&...&0\\
	\vdots&\vdots&\ddots&\vdots\\
	mk\sqrt{2}A(G)&0&...&0
\end{bmatrix}_{(mn+n)}$\\

\hspace{0.72in}$=\begin{bmatrix}
	mk\sqrt{2}&mk\sqrt{2}&...&mk\sqrt{2}\\
	mk\sqrt{2}&0&...&0\\
	\vdots&\vdots&\ddots&\vdots\\
	mk\sqrt{2}&0&...&0
\end{bmatrix}_{(m+1)} \otimes A(G)$\\

Let $C=\begin{bmatrix}
	mk\sqrt{2}&mk\sqrt{2}&...&mk\sqrt{2}\\
	mk\sqrt{2}&0&...&0\\
	\vdots&\vdots&\ddots&\vdots\\
	mk\sqrt{2}&0&...&0
\end{bmatrix}_{(m+1)}$ and let eigenvalues of matrix $C$ be $\delta_1, \delta_2, ..., \delta_{(m+1)}.$ Since $B$ is of rank two therefore there two non zero eigenvalues $\delta_1$ and $\delta_2$ (say). Then

\begin{center}
		$ tr(B)=\delta_1+\delta_2 $
\begin{equation}
\delta_1+\delta_2=mk\sqrt{2}
\end{equation}
	$\delta_1+\delta_2=mk\sqrt{2}$\\
	and	
	$  tr(B^2)=\delta_1^2+\delta_2^2 $
	\begin{equation} \delta_1^2+\delta_2^2=4m^3k^2+2m^2k^2 \end{equation} 
\end{center}
Solving $(3)$ and $ (4) $ we get\\
\begin{center}
	$\delta^2-mk\sqrt{2}\delta -2m^3k^2=0$ 
\end{center}	
	whose roots are $ \delta_1 $ and $ \delta_2 $

\noindent  $\therefore ~~~~$ Characteristic equation of $B$ is
\begin{center}
	$Char(B)=\delta^{(m-1)}(\delta^2-mk\sqrt{2}\delta-2m^3k^2)=0$.
\end{center}
Therefore,
\begin{center}

	$\implies \delta=0$ and $\displaystyle\delta=mk\sqrt{2}(\frac{1\pm \sqrt{4m+1}}{2})$.
\end{center}
\noindent Therefore,
\begin{center}
	$\displaystyle spec(C)=\begin{array}{ccc}
		\left(
		\begin{array}{ccc}
			0&mk\sqrt{2}(\frac{1+\sqrt{4m+1}}{2})&mk\sqrt{2}(\frac{1- \sqrt{4m+1}}{2})\\
			m-1&1&1
		\end{array}
		\right)
	\end{array}$
\end{center}
\noindent by preposition [2]
\begin{center}
	$\displaystyle spec (S(D_m(G)))=\begin{array}{ccc}
		\left(
		\begin{array}{ccccccccc}
			0\lambda_1&...&0\lambda_n&p\lambda_1&...&p\lambda_n&q\lambda_1&...&q\lambda_n\\
			m-1&...&m-1&1&...&1&1&...&1
		\end{array}
		\right)
	\end{array}$
\end{center}
Where, $\displaystyle p=mk\sqrt{2}(\frac{1+\sqrt{4m+1}}{2})$ and $\displaystyle q=mk\sqrt{2}(\frac{1- \sqrt{4m+1}}{2})$.\\

\noindent Therefore,

\begin{center}
	$\displaystyle ES(D_m(G))=\sum\limits_{i}^{n}\left|mk\sqrt{2}(\frac{1\pm \sqrt{4m+1}}{2})\lambda_i
	\right|$\\
	
\hspace{1.9in}	$\displaystyle =mk\sqrt{2}\sum\limits_{i}^{n}\left|\lambda_i\right|\left[\frac{\sqrt{4m+1}+1}{2}+\frac{\sqrt{4m+1}-1}{2}\right]$\\

\hspace{0.1in}	$\displaystyle = mk\sqrt{8m+2}\sum\limits_{i}^{n}\left|\lambda_i\right|$\\
$\displaystyle = mk\sqrt{8m+2} \varepsilon(G)$
\end{center}

	\begin{table}[h!]
	\caption{Sombor index and Sombor energy of $m$-splitting and $m$-shadow graph of some specific regular graph.}
	\begin{center}
	\scalebox{0.9}{

		\begin{tabular}{|c|c|c|c|c|c|}
		\hline
	Graph & $SO(G)$ & $SO(Spl_m(G))$ & $SO(D_m(G))$  & $SE(Spl_m(G))$ & $SE(D_m(G))$  \\
	$(G)$&&&&&\\
	\hline
	$C_n$&$2\sqrt{2}n$&$2\sqrt{2}n(m\sqrt{2m^2+4m+4}+(m+1))$&$2\sqrt{2}n(m^3+m^2)$&$2\sqrt{2}(m+1) \varepsilon(C_n)$&$2m\sqrt{8m+2} \varepsilon(C_n)$\\
	\hline
	$K_n$&$\frac{n(n-1)^2}{\sqrt{2}}$&$\frac{n(n-1)^2}{\sqrt{2}}(m\sqrt{2m^2+4m+4}+(m+1))$&$\frac{n(n-1)^2}{2}(m^3+m^2)$&$2\sqrt{2}(n-1)^2(m+1)$&$2(n-1)^2 m\sqrt{8m+2}$\\
	\hline
	$Q_n$&$2^{(n-\frac{1}{2})}n^{2}$&$2^{(n-\frac{1}{2})}n^{2}(m\sqrt{2m^2+4m+4}+(m+1))$&$2^{(n-\frac{1}{2})}n^{2}(m^3+m^2)$&$2\sqrt{2}n(m+1)\lceil \frac{n}{n} \rceil{n\choose {\lceil \frac{n}{n} \rceil}}$&$2mn\sqrt{8m+2}\lceil \frac{n}{n} \rceil{n\choose {\lceil \frac{n}{n} \rceil}}$\\
	\hline
	$K_{n,n}$&$\sqrt{2}n^3$&$\sqrt{2}n^3(m\sqrt{2m^2+4m+4}+(m+1))$&$\sqrt{2}n^3(m^3+m^2)$&$2\sqrt{2}n^2(m+1)$&$2mn^2\sqrt{8m+2}$\\
%	Petersen Graph&&&&&\\

	\hline
	
	\end{tabular}}
\end{center}
	
\end{table}

\noindent\textbf{Conclusion.}
Recently I. Gutman introduced new degree based topological index called as Sombor index. We have obtained general formula for Sombor index and relation between energy and Sombor energy of $m$-splitting and $m$-shadow graphs of $k$-regular graph.

\end{document}